\documentclass[georeport,12pt]{geophysics}

\usepackage{natbib}
\usepackage{graphicx}
\usepackage{color}
\usepackage{amsmath}
\usepackage{amsthm}
\usepackage{amssymb}
\usepackage{amsbsy}
\usepackage{wasysym}
\usepackage{setspace}
\usepackage{bm}
\usepackage{float}
\usepackage{hyperref}
\usepackage{algorithm}
\usepackage{algpseudocode}

\newcommand{\bx}{\mathbf{x}}

\newcommand{\bv}{\mathbf{v}}

\newcommand{\bR}{\mathbf{R}}

\newcommand{\Ja}{J_{\alpha}}
\newcommand{\tJa}{\tilde{J}_{\alpha}}
\newcommand{\tJz}{\tilde{J}_{0}}
\newcommand{\Jase}{J_{\rm MSWI,\alpha,\sigma,\epsilon}}
\newcommand{\tJase}{\tilde{J}_{\rm MSWI,\alpha,\sigma,\epsilon}}
\newcommand{\JAase}{J_{\rm AWI,\alpha,\sigma,\epsilon}}
\newcommand{\tJAase}{\tilde{J}_{\rm AWI,\alpha,\sigma,\epsilon}}

\newcommand{\ua}{u_{\alpha}}
\newcommand{\uase}{u_{\rm MSWI,\alpha,\sigma,\epsilon}}
\newcommand{\uaase}{u_{\rm AWI,\alpha,\sigma,\epsilon}}
\newcommand{\va}{v_{\alpha}}
\newcommand{\uz}{u_0}
\newcommand{\wl}{w_{\lambda}}
\newcommand{\fl}{f_{\lambda}}
\newcommand{\Fl}{F_{\lambda}}
\newcommand{\Fla}{F^0_{\lambda}}
\newcommand{\Flb}{F^1_{\lambda}}
\newcommand{\Sl}{S_{\lambda}}
\newcommand{\Sla}{S^0_{\lambda}}
\newcommand{\Slb}{S^1_{\lambda}}
\newcommand{\dl}{d_{\lambda}}
\newcommand{\dzl}{d^0_{\lambda}}
\newcommand{\dla}{d^0_{\lambda}}

\newcommand{\fwl}{\hat{w}_{\lambda}}
\newcommand{\ul}{u_{\sigma}}
\newcommand{\uzl}{u^0_{\sigma}}
\newcommand{\Te}{T_{\epsilon}}
\newcommand{\Tekds}{T_{\epsilon,\kappa,d,\sigma}}
\newcommand{\fu}{\hat{u}}

\newtheorem{lemma}{Lemma}
\newtheorem{theorem}{Theorem}
\newtheorem{prop}{Proposition}

\newtheorem{cor}{Corollary}

\setfigdir{Fig}

\begin{document}
\title{Adaptive Waveform Inversion for Transmitted Wave Data}
\author{William W. Symes}

\begin{abstract} Adaptive Waveform Inversion (AWI) applied to
  transient transmitted wave data can yield estimates of index of
  refraction (or wave velocity) similar to those obtained by travel
  time inversion. The AWI objective function measures normalized
  mean-square dispersion about zero time of a family of filters, one
  filter for each source-reeciver pair, designed to match predicted
  data to observed data. Provided that the data contain a single
  smooth wavefront, this function approaches the mean square
  traveltime error as data wavelength tends to zero. The scaling of
  each filter to have unit norm is responsible for the simple relation
  with travel time tomography. A similar approach, Matched Source
  Waveform Inversion (MSWI), does not normalize the filters and has a
  looser relation with mean-square travel time error. If substantial
  energy is spread amongst multiple arrivals, on the other hand,
  neither AWI nor MSWI objectives approximate the travel time mean
  square erorr, and attempts to minimize them by local
  (Newton-related) optimization are as likely to stagnate at models
  predicting erroneous travel times as is least-square Full Waveform
  Inversion (FWI). The matching filter at the heart of this approach
  must be ``pre-whitened'', that is, computed by solving a regularized
  least squares problem for filtered data misfit. In order for the
  relation to traveltime tomography to hold, and for the filtered data
  misfit to stay relatively small, the regularization weight must be
  coupled to data wavelength. These objective functions can also be
  used as penalty terms in a homotopy of objectives connecting AWI or
  MSWI with FWI. \end{abstract}

\section{Introduction}
Full waveform inversion (FWI) refers to the inference of mechanical
parameter distribution within a material object  (the earth, a human
body, a manufactured artifact,...) from remote measurements of waves
propagating through it. These mechanical parameters appear as
coefficients in systems of partial differential equations governing
wave propagation, and numerical solutions of these wave equations can
be used to predict experimental data. FWI consists in
updating the parameters to minimize a numerical measure of  the difference between computationally
predicted wave data and actual observed
data. Such data is typically oscillatory. If wave velocities or their parametric equivalents are amongst the parameters
to be estimated, then updates cause shifts in the arrival times of
data oscillations, which in turn causes commonly used numerical error measures
to oscillate. This ``cycle skipping'' phenomenon may cause iterative local optimization algorithms (relatives of
Newton;'s method, the only practical algorithmic choices) to stagnate
far from good data fit, unless the initial
parameter estimates accurately predict the arrival times of waves \cite[]{GauTarVir:86,VirieuxOperto:09,HuetalIEEE:18}.

Adaptive Waveform
Inversion \cite[]{Warner:16} (AWI) is a modification of FWI, designed to overcome the
cycle-skipping pathology. AWI has been successfully applied to
industry-scale seismic data
\cite[]{GuaschWarnerRavaut:GEO19,Warneretal:SEG21} and shows promise
for ultrasonic imaging of the the human brain
\cite[]{Guaschetal:NPJDM20}. These and other works discuss heuristic justifications for AWI
and many details of implementation. However, a clear explanation for its
effectiveness is so far missing from the literature.

This paper supplies such an explanation for a special case of AWI, the
recovery of smooth acoustic medium parameters from waves generated at
localized energy sources, propagated through a fluid, and recorded at
receiver locations remote from the sources (transmission data). I show that the AWI
objective function is an asymptotic (in time frequency) approximation
to the mean-square time-of-travel misfit, provided that arrival times
of transmitted waves are single. That is, a unique ray of geometric
asymptotics connects each sources-receiver pair, as is the case if the associated
Riemannian metric is simple in a domain containing the sources and
receivers \cite[]{PestovUhlmann:05,StefanovUhlmann:05}. The error in
this approximation is proportional to the source mean-square
wavelength. I also show that this relation fails if arrival times have
higher multiplicity. Travel time tomography (minimization of
mean-square travel time error) appears to be much less sensitive to
choice of initial model than is FWI. In fact, it is commonly used to
construct an initial model for FWI
\cite[]{Bordingetal:87,SirguePratt:04,VirieuxOperto:09}. At least for
acoustic transmission data with single arrivals, AWI implicitly
combines travel time tomography with FWI.

The arguments presented here for AWI are refinements of those given by
\cite{HuangSymes:Geo17} for Matched Source Waveform Inversion (MSWI),
a similar modification of FWI
\cite[]{HuangSymes2015SEG,HuangSymes:Geo17}. which were in turn based
on ideas presented in Song's thesis \cite[]{Song:94}. Since it
involves no additional effort, I state results for MSWI in parallel
with those for AWI.

The next section establishes the class of problems to be studied and notations
for their description, defines the AWI and MSWI objective functions,
and states the main results of this paper. In the following section, I present the
asymptotic analysis leading to the main results. In
effect, AWI regularizes or relaxes the FWI problem. For many inverse
problems, regularization has been accomplished via additive penalty
functions. In section four, I explain the relation between AWI, MSWI, and 
penalty formulations of the wave inverse problem. The close relation
between AWI objective and travel time misfit is limited in scope: in
the fifth section, I show that it does not hold in general for
transmitted wave data in which multiple rays connect source and
receiver (multiiple wave arrivals, non-simple metric).

The analysis explained here sheds no light on the effectiveness of AWI
for reflected wave data, and the approach has so far been limited to
inversion of data modeled by acoustic wave propagation. I discuss
these and other aspects of AWI and related approaches in a brief final
section.

\section{Main Results}

The version of AWI introduced by \cite{Warner:16} is based on an
active source acoustic model of data generation in Euclidean space
$\bR^3$. 

Each acoustic source is modeled an isotropic point radiator
with known location $\bx_s\in \bR^3$. These
sources generate acoustic waves traveling through Euclidean space. The
data are samples of these waves at receivers, modeled as isotropic
point sensors, recording pressure at points
${\bf x}_r \in \bR^3$. The collection $X_{sr}$ of source-receiver pairs
$(\bx_s,\bx_r)$ for which data is recorded should not intersect the
diagonal - that is, no receiver should be co-located with a source for
which it is active. The set of source locations $X_s = \{\bx_s: (\bx_s,\bx_r) \in X_{sr}
\mbox{ for some } \bx_r \in \bR^3\}$ will also be useful.

Throughtout the following discussion, it is assumed that
$X_{sr}$ is finite. 
  
The pressure and velocity fields $p({\bf x},t;{\bf x}_s)$,
${\bf v}({\bf x},t;{\bf x}_s)$ for the source location ${\bf x}_s$
depend on the bulk modulus $\kappa$, buoyancy $\beta$
(reciprocal of the density $\rho$), and source wave
form $f \in C_0^{\infty}(\bR)$ through the acoustic system
\begin{eqnarray}
  \label{eqn:awe}
 \frac{\partial p}{\partial t} & = &- \kappa \nabla \cdot {\bf v} +
                                    f(t) \delta({\bf x}-{\bf x}_s); \nonumber \\
\frac{\partial {\bf v}}{\partial t} & = & - \beta \nabla p, \nonumber \\ 
p, {\bf v} & = & 0 \mbox{ for }  t \ll 0.
\end{eqnarray}
For simplicity, I assume that the system \ref{eqn:awe} holds over all of Euclidean
space-time $(\bx,t) \in \bR^3 \times \bR$, and ignore physical
boundaries such as the Earth's surface. Likewise for simplicity, I
assume in this study that the material density $\rho$ is spatially
homogeneous and known.

The parameter to be determined in the inverse problem studied here is
the bulk modulus $\kappa$. It is presumed is bounded and smooth: $\log
\kappa \in M$, where $M$ is a bounded subset of $C^{\infty}(\bR^3) \cap
L^{\infty}(\bR^3)$, on which further constraints will be placed by and
by. 

Under these hypotheses, the symmetric hyperbolic system \ref{eqn:awe}
has a unique distribution solution (for example \cite{Lax:PDENotes},
Theorem 7.4). The Propagation of
Singularities Theorem (\cite{Tay:81}, Theorem VI.2.1) 
imples that $p(\cdot,\cdot;\bx_s)$ is smooth away from $\{\bx_s\}
\times \bR$, so that the traces $p(\bx_r, \cdot;\bx_s)$ on
$\{(\bx_r,t): t \in \bR\}$ are well-defined.

Define the {\em modeling operator}
\begin{equation}
  \label{eqn:modop}
  F: exp(M)  \rightarrow \Pi_{X_{sr}} C^{\infty}(\bR)
\end{equation}
by
\begin{equation}
  \label{eqn:modopdef}
  F[\kappa]= p|_{X_{sr} \times \bR}
\end{equation}

In this setting, acoustic FWI means: given the set of
source-receiver position pairs $X_{sr}$, source wavelet
$f \in C_0^{\infty}(\bR)$,  density $\rho$, and data traces
$d \in \Pi_{X_{sr}} C^{\infty}(\bR)$, find a bulk modulus $\log
\kappa$
(or equivalently wave velocity $c=\sqrt{\kappa/\rho}$)
so that $F[\kappa] \approx d$ in some sense. If
$d(\bx_r,\cdot;\bx_s)$ and $F[\kappa](\bx_r,\cdot;\bx_s)$ lie in a normed subspace of ${\cal D}'(\bR)$ for each
$(\bx_s,\bx_r) \in X_{sr}$, then the
norm provides a straightforward candidate for this sense of
approximation.

The analysis to follow requires three further constraints on the set
of source-receiver location pairs $X_{sr}$ and the domain
$M \subset C^{\infty}(\bR^3) \cap L^{\infty}(\bR^3)$ of $F$. The first is
the {\em Exponential Decay} assumption:
\begin{quote}
  For any $n > 0$, $f \in C_0^{\infty}(\bR)$,
  there exist $K >
  0, T \in \bR,$ and $\delta>0$ so
  that for any $k$ with $0 \le k \le n$, $(\bx_s,\bx_r) \in
 X_{sr},$ and any $\log \kappa \in M$, the solution $(p,\bv)$ of \ref{eqn:awe} is smooth
  for $t>T$ and satisfies
  \[
    |\partial_t^k p(\bx_r,t;\bx_s)| \le K e^{-\delta t}, \, t > T.
  \]
\end{quote}
SInce $X_{sr}$ is bounded, the property described in this assumption
is a special case of {\em local energy decay}, a much studied topic. The {\em
  non-trapping} hypothesis is sufficient to ensure local energy decay:
all rays of geometric asymptotics (that is, geodesics for the Riemannian metric associated
to the system \ref{eqn:awe}) escape to infinity (\cite{Hristova:09},
Theorem 2, \cite{EgorovShubin}, Theorem 2.104, p. 192). In particular,
the Exponential Decay assumption is satisfied if for any $R > r>0$ and $\log
\kappa \in M$, there exists $T>0$ so that any ray ${\bf X}:\bR^+\rightarrow \bR^3$
  with ${\bf X}(0) \in B_r(0)$ satisfies $|{\bf X}(t)|>R$ for $t > T$. That is,
exponential decay may follow from a propery of $M$.

The second and third assumptions explicitly involve constraints on ray geometry
(therefore on $\kappa$). A {\em normal neighborhood} of
$\bx_s \in \bR^3$ is an open subset $\Omega$ of $\bR^3$ 
within which every point is connected to $\bx_s$ by  a unique ray  (\cite{Friedlander:75},
p. 16). Therefore, to every point $\bx \in \Omega$ is associated a
unique time of travel (geodesic distance) $\tau(\bx_s,\bx)$ from $\bx_s$ to $\bx$.
The {\em Single Arrival} assumption is
\begin{quote}
  $\Omega\subset \bR^3$ and $M \subset C^{\infty}(\bR^3) \cap
  L^{\infty}(\bR^3)$ are bounded open sets for which $X_{sr} \subset
  \Omega \times \Omega$, and for every source point $\bx_s\in X_{s}$
  and $\log \kappa \in M$, $\Omega$ is a normal neighborhood of
  $\bx_s$.
\end{quote}

The {\em No Return} assumption is
\begin{quote}
  for every source point $\bx_s \in X_{s}$, $\log \kappa \in M$, ray ${\bf X}:\bR^+\rightarrow \bR^3$
  with ${\bf X}(0)=\bx_s$ and any $t>0$, either $\{{\bf X}(s): 0\le s \le t\} \subset \Omega$
  or $\{{\bf X}(s): t < s < \infty\} \subset \bR^3 \setminus \Omega$.
\end{quote}
This assumption states that once a ray leaves
$\Omega$, its subsequent intersection with $\Omega$ is null (hence the name).

Thanks to Exponential Decay, each trace $p(\bx_r,\cdot;\bx_s)$ is in $L^{p}(\bR)$ for $p
\ge 1$.
Denote by $L^p(X_{sr} \times \bR)$ the product space $l^p(X_{sr},L^p(\bR))$
with the usuall product norm. 

Accordingly, FWI may be posed as a least-squares
problem: given data $d \in
L^2(X_{sr}\times \bR)$, find $\kappa$ to minimize
\begin{equation}
  \label{eqn:fwi}
  J_{\rm FWI}[\kappa;d]= \frac{1}{2}\|F[\kappa]-d\|^2_{L^2(X_{sr} \times \bR)},
\end{equation}
As mentioned earlier, optimization algorithms related to Newton's method, applied to
$J_{\rm FWI}$, tend to stagnate far from global minimizers if
the initial choice of $\kappa$ does not predict the arrival times of
waves correctly within a half-wavelength or so (the meaning of
``wavelength'' in this context will be explained below)
\cite[]{GauTarVir:86,VirieuxOperto:09}.

AWI modifies the least-squares FWI problem by interposing a convolution operator or
{\em filter} $u \in L^2(X_{sr} \times \bR)$. That is, the residual
data $F[\kappa]-d$ is replaced by the filtered residual data $F[\kappa]*u-d$.
Since $F[\kappa] \in L^1(X_{sr}\times \bR)$, the filtererd residual
makes sense.
The filter $u$ is intended to shift localized wave energy in $F[m]$ to
match that present in $d$. If on the other hand $u(\bx_r,t;\bx_s) =
\delta(t)$ for each $(\bx_x,\bx_r) \in X_{sr}$, then $tu = 0$. This
observation motivates the definition of the AWI objective function
\begin{equation}
  \label{eqn:awidefmain}
  J_{\rm AWI,\sigma}[\kappa;d] = \sum_{(\bx_s,\bx_r) \in X_{sr}}
  \frac{\|Tu_{\sigma}[\kappa;d](\bx_r,\cdot;\bx_s)\|^2}{\|u_{\sigma}[\kappa;d](\bx_r,\cdot;\bx_s)\|^2}.
\end{equation}
in which $Tu(\bx_r,t;\bx_s) = tu(\bx_r,t;\bx_s)$ and
\begin{equation}
  \label{eqn:reg}
  u_{\sigma}[\kappa;d] = \mbox{argmin}_{u \in L^2(X_{sr} \times \bR)} (\|F[\kappa]*u-d\|^2 + \sigma\|u\|^2)
\end{equation}
The objective function for the closely related MSWI approach,
mentioned in the introduction, is  
\begin{equation}
  \label{eqn:mswidefmain}
  J_{\rm MSWI,\sigma}[\kappa;d] = \sum_{(\bx_s,\bx_r) \in X_{sr}}
  \|Tu_{\sigma}[\kappa;d](\bx_r,\cdot;\bx_s)\|^2.
\end{equation}
Note that while the Lax-Milgram Theorem (\cite{Yosida}, section
III.7) ensures the existence of a unique solution
$u_{\sigma}[\kappa;d] \in L^2(X_{sr}\times \bR)$ of the least squares
problem \ref{eqn:reg},
it is not immediiately obvious that $Tu_{\sigma}[\kappa;d] \in L^2(X_{sr}\times \bR)$.

Asymptotic analysis of AWI and MSWI is based on a family of source
wavelets: for $\lambda > 0$,
\begin{equation}
  \fl(t)= \lambda^{\frac{1}{2}}f_1\left(\frac{t}{\lambda}\right).
\end{equation}
The relation between the dimensionless scale parameter $\lambda$ and
RMS wavelength will be explained below: $\lambda$ may be regarded as a
proxy for wavelength. Correspondingly, denote the
modeling operator with $f=\fl$ by $\Fl$.

The main result of this analysis is

\begin{theorem}
  \label{thm:main}
  Assume that Exponential Decay, Single Arrival, and No Return
  hold for $X_{sr} \subset \Omega \times \Omega, M \subset C^{\infty}(\bR^3) \cap L^{\infty}(\bR^3)$.
  Select $\log \kappa^* \in M$. For $0 <
  \lambda \le 1$, set
  $\dl=\Fl[\kappa^*]$ (``noise-free data''), and  define $\ul[\kappa,\dl]$
  to be the solution of 
  equation \ref{eqn:reg} with $F=\Fl$ and $d=\dl$..
  
  \begin{itemize}
  \item[1. ] for $(\bx_s,\bx_r) \in X_{sr}$, $\log \kappa \in M$, and
    $\sigma, \lambda >0$, $Tu_{\sigma}[\kappa;\dl](\bx_r,\cdot;\bx_s) \in
    L^2(\bR)$;
  \item[2. ] There exist $R_*,R^*,\lambda_*$ with $0<R_*\le R^*, 0 <
    \lambda_*\le 1$ depending on
    $M, X_{sr}$ so that
    \begin{itemize}
    \item[a. ] $\|\Sl[\kappa]\ul[\kappa;\dl]-\dl\| \le
      \frac{1}{2}\|\dl\|$ implies that $\sigma \le R^*\lambda$;
    \item[b. ] $\sigma \le R_*\lambda$ and $\lambda \le \lambda_*$
      implies that $\|\Sl[\kappa]\ul[\kappa;\dl]-\dl\| \le
      \frac{1}{2}\|\dl\|$.
    \end{itemize}
  \end{itemize}
  Choose $r > 0$ and assume that
  \begin{equation}
    \label{eqn:lamsig}
    \sigma = r \lambda.
  \end{equation}.
  Then there exist $0 < C_* \le C^*$, $0 < C$, and $W: M \times X_{sr}
  \rightarrow \bR^+$ with $C_* \le
  W \le C^*$ so that for
  $\log \kappa \in M$, $(\bx_s,\bx_r) \in X_{sr}$,
  \begin{itemize}
  \item[3. ] for $0 <  \lambda \le 1$, 
    \begin{equation}
      \label{eqn:mswitt}
      |\|T\ul[\kappa;\dl](\bx_r,\cdot;\bx_s)\|- \lambda^{-\frac{1}{2}} W[\kappa^*,\kappa]^{\frac{1}{2}}(\bx_s,\bx_r)|\tau[\kappa](\bx_s,\bx_r)-\tau[\kappa^*](\bx_s,\bx_r)|| \le C
      \lambda^{\frac{1}{2}};
    \end{equation}
  \item[4. ]for sufficiently small $\lambda$, 
    \begin{equation}
      \label{eqn:awitt}
      \left|\frac{\|T\ul[\kappa;\dl](\bx_r,\cdot;\bx_s)\|}{\|\ul[\kappa,\dl](\bx_r,\cdot;\bx_s)\|}-
        |\tau[\kappa](\bx_s,\bx_r)-\tau[\kappa^*](\bx_s,\bx_r)|\right| \le C \lambda
    \end{equation}
  \end{itemize}
\end{theorem}
\noindent{\bf Remark:} The heuristic justification for AWI is the possibility of keeping the
filtered predicted data $\Fl[\kappa]*\ul[\kappa;d]$ close to the data
$d$ while updating the bulk modulus $\kappa$. The conclusions 2a and
2b show that for this to be the case, the regularization parameter $\sigma$
must be dominated by the wavelength proxy $\lambda$. This observation underlies the choice
$\sigma = r\lambda$ used in conclusions 3 and 4. A looser relation
would clearly suffice, however the argument lis slightly cleaner if
$\sigma$ is simply held proportional to $\lambda$ as $\lambda
\rightarrow 0$. Note that the units of $r$ and $\sigma$ are pressure$^2 \times$time$^2$.

Summation over the source-receiver positions in $X_{sr}$ yields
\begin{cor}
  \label{thm:cor1}
  In addition to the assumptions of Theorem \ref{thm:main}, define $J_{\rm
    AWI,\sigma}$ and $J_{\rm MSWI,\sigma}$ by equations
  \ref{eqn:awidefmain}, \ref{eqn:mswidefmain} respectively. Then
\begin{equation}
\label{eqn:cor1mswi}   
\lim_{\lambda \rightarrow 0} \lambda J_{\rm MSWI,r\lambda}[\kappa;\dl]
= \sum_{(\bx_s,\bx_r) \in X_{sr}} W[\kappa^*,\kappa](\bx_s,\bx_r) |\tau[\kappa](\bx_s,\bx_r) - \tau[\kappa^*](\bx_s,\bx_r)|^2.
\end{equation}
and
\begin{equation}
\label{eqn:cor1awi}   
\lim_{\lambda \rightarrow 0} J_{\rm AWI,r\lambda}[\kappa;\dl]
= \sum_{(\bx_s,\bx_r) \in X_{sr}}  |\tau[\kappa](\bx_s,\bx_r) - \tau[\kappa^*](\bx_s,\bx_r)|^2.
\end{equation}
\end{cor}

\section{Analysis of AWI and MSWI}
This section is devoted to the proofs of the main result, Theorem
\ref{thm:main}.  The principal tool in the analysis is Hadamard's
decomposition of the acoustic Green's function, described in the
second
subsection below, along with some of its consequences, the chief being the
square-integrability of the time-scaled regularized filter
$tu_{\sigma}$. The next subsection reveals the link to travel time via
an approximate solution of the regularized least squares problem
\ref{eqn:reg}, based on the leading term in Hadamard[s decomposition.
The final subsection estimates the remainder and finishes the proof of
the main result.

\subsection{A Note on Notations  for Norms and Bounds}
The Banach spaces
$l^p(X_{sr},L^p(\bR))$ will be abbreviated as $L^p(X_{sr} \times
\bR)$. The norm in a Banach space $B$ will be denoted
$\|\cdot\|_{B}$. The norm symbol
$\|\cdot\|$ without subscript will be understood to mean either
$\|\cdot\|_{L^2(\bR)}$ or $\|\cdot\|_{L^2(X_{sr}\times
  \bR)}$, as context indicates.

In a statement pertaining to all
$\kappa \in M$, the letter $C$, or sometimes $C_k,
k=0,1,2,...$ or
$C_*,C^*$, will designate a constant (or constants) independent of
$\kappa$.

Generally, constants denoted $C$ etc. will absorb any dimensional
scale factors necessary to make sense of the expressions in which they
appear. I will call out those few cases in which explicit use of
dimensinoal scale factors is necessary to avoid irretrievable
dimensional confusion.

\subsection{Hadamard's Decomposition}
The Hadamard decomposition of the acoustic Green's function \cite[]{Friedlander:75,Qian:JCP24} provides a
basis for analysis of AWI and MSWI. As presented by
\cite{Friedlander:75}, for instance, this construction applies to
second-order wave equations, but the system \ref{eqn:awe} is
equivalent fo a second-order problem: differentiate the first equation
with respect to $t$, eliminate $\bv$ via the second equation, and
rearrange to obtain
\begin{equation}
  \label{eqn:awe2}
  \frac{\rho}{\kappa}\frac{\partial^2 p}{\partial t^2} - \nabla^2 p -
  \nabla \log \rho \cdot \nabla p = \rho \frac{\partial f}{\partial t}
  \delta(\bx-\bx_s), \nonumber \\
  p=0, \, t \ll 0.
\end{equation}
If $(p,\bv)$ is a distribution solutions of the system
\ref{eqn:awe} then $p$ is a solution of equation
\ref{eqn:awe2}, as follows by differentiating the first equation of
\ref{eqn:awe} with respect to $t$ and eliminating $\bv$ via the second
equation. Conversely, if $p$ is a solution of the second-order initial
value problem\ref{eqn:awe2}, then under the conditions imposed in the
previous section,
\[
  \bv = \beta \int_{-\infty}^t p
\]
defines a distribution, and $(p,\bv)$ is a weak solutino of the system \ref{eqn:awe}.

Hadamard's decomposition is local, in both space-time and in the set of
acoustic models. In the present context, it requires the Single
Arrival hypothesis, and for global time validity the No Return
assumption. 

The construction described by \cite{Friedlander:75} produces for each $\log \kappa \in M$ functions $\tau[\kappa]$ (travel time), $a[\kappa]$
(geometric amplitude) on $X_s \times \Omega$ and $b[\kappa]$
(remainder kernel) on $X_s \times \Omega \times \bR$, smooth off the
diagonal - that is, for each $\bx_s \in X_s$, $\tau[\kappa](\bx_s,\cdot) \in C^{\infty}(\Omega \setminus
  \{\bx_s\})$, and similarly for $a[\kappa],b[\kappa]$; The traveltime $\tau$
  satisfies the eikonal equation
  \[
    \sqrt{\frac{\kappa}{\rho}} |\nabla \tau(\bx_s,\cdot)| = 1 \mbox{ in } \Omega \setminus
    \{\bx_s\}, \, \lim_{\bx \rightarrow \bx_s} \tau(\bx_s,\bx) = 0.
  \]
A by-product of the construction is a bound $A>0$ so that for any $\log \kappa \in M$
and $(\bx_s,\bx_r) \in X_{sr}$,
\begin{equation}
  \label{eqn:aest}
  |\log a[\kappa](\bx_s,\bx_r)| <  A.
\end{equation}
With no further assumptionss, there is also a uniform bound on
$b(\bx_r,t;\bx_s)$ and
its derivatives over any finite $t$-interval. However, with the
addition of the Exponential Decay assumption, the bound becomes
global in time: for any $k \ge 0$, there exists $B \ge 0,\delta > 0$
so that for any $(\bx_s,\bx_r) \in X_{sr}$, $\log \kappa \in M$,
\begin{equation}
  |\partial^k_t b[\kappa](\bx_r,t;\bx_r,s)| \le B e^{-\delta t}
  \label{eqn:best}
\end{equation}

The Hadamard
decomposition of the Green's function $G[\kappa]$, i.e. the solution of the
equation \ref{eqn:awe2} for $f = \beta H$, is::
\begin{equation}
  \label{eqn:green}
  G[\kappa](\bx_s,\bx,t) = a[\kappa](\bx_s,\bx)\delta(t - \tau[\kappa](\bx_s,\bx)) +
  b[\kappa](\bx_s,\bx,t-\tau[\kappa])H(t - \tau[\kappa](\bx_s,\bx))
\end{equation}
Define $\delta_{\tau} = \delta(t-\tau)$ to be the shifted delta, that
is, $\langle \delta_{\tau},\phi\rangle = \phi(\tau)$. A useful
restatement of the definition \ref{eqn:green}, suppressing the
space-time coordinates, is
\begin{equation}
  \label{eqn:regreen}
  G[\kappa] = (a[\kappa]\delta +  (b[\kappa]H))*\delta_{\tau[\kappa]}.
\end{equation}

\subsection{Tihonov Regularization}
Introduce the notation $S[\kappa]$ for the operator of convolution
with $F[\kappa]$, that is,
\begin{equation}
  \label{eqn:sdef}
  S[\kappa]u = u*F[\kappa]
\end{equation}
$S$ represents an
extension of $F$, in the sense that with the choice
$u(\bx_r,t;\bx_s)=\delta(t)$, the output of $F$ is recovered:
\begin{equation}
  \label{eqn:sconsist}
  S[\kappa]\delta = F[\kappa].
\end{equation}
The regularized least-squares problem \ref{eqn:reg} can be re-written
\begin{equation}
  \label{eqn:regs}
    u_{\sigma}[\kappa;d] = \mbox{argmin}_{u \in L^2(X_{sr} \times \bR)} (\|S[\kappa]u-d\|^2 + \sigma\|u\|^2)
  \end{equation}
\cite{Warner:16} present the definition of the AWI objective function
\ref{eqn:awidefmain} initially as if $\sigma=0$, in fact presuming a
zero-residual solution (that is, $S[\kappa]u=d$) but later note that
``pre-whitening'' is necessary for a well-behaved filter construction.
``Pre-whitening'' is a synonym for Tihonov regularization, that is,
the least-squares minimization \ref{eqn:regs} with $\sigma>0$.

Regularization has another benefit, in that it makes the objective
well-defined. Indeed, if the data is noise-free, that is,
$d=F[\kappa]=S[\kappa]\delta$, then $u_0[\kappa,d]=\delta$, but the
denominator in the definition \ref{eqn:awidefmain} is not defined.
Of course, the numerator vanishes, motivating the definition of AWI,
but the objective function is nonsense.

For future reference, examine the implications of definition
\ref{eqn:regs} for the units of the regularization weight $\sigma$,
which is necessarily dimensional. Since the data $d$ has dimensions of
pressure, so must $S[m]u = F[m]*u$. Of course $F[m]$ is also sampled
pressure, and convolution is a $t$ integration, so for $S[m]u$ to be a
pressure, $u$ must have dimensions of 1/time. Therefore
$\|S[m]u-d\|^2$ has dimensions of pressure$^2 \times $time, $\|u\|^2$
has dimensions of 1/time, whence $\sigma$ must be dimensioned as
pressure$^2 \times$ time$^2$.

From \ref{eqn:green} and the
definition \ref{eqn:sdef} of  $S[\kappa]$,
\[
  S[\kappa]u = \rho( a[\kappa]\delta    +  (b[\kappa]H)) *\delta_{\tau[\kappa]}*\partial_t f *u,
\]
\begin{equation}
  \label{eqn:sgood}
  = \rho(a[\kappa]\delta_{\tau[\kappa]}*\partial_t f +
  (b[\kappa](0)\delta + (\partial_t b[\kappa])H)
  *\delta_{\tau[\kappa]}*f) * u.
\end{equation}
The formal adjoint $S[\kappa]^T$ is given by
\begin{equation}
  \label{eqn:sadjgood}
  S[\kappa]^Tu = \rho( a[\kappa]\delta    -  (\check{b}[\kappa]\check{H}))
 *\delta_{-\tau[\kappa]}*\partial_t \check{f} *u,
\end{equation}
in which $\check{u}(t) = u(-t)$.
The key properties of $S$ follow directly from Young's inequality and
the representation \ref{eqn:sgood}:
\begin{lemma}
  \label{thm:snorm}
  Under the Exponential Decay, Single Arrival, and No Return assumptions, the
  expression \ref{eqn:sgood} define
  bounded operators on
  $\Pi_{X_{sr}} {\cal S}(\bR)$ and on
  $L^2(X_{sr} \times \bR)$, and there exists $C>0$ so
  that for any $\kappa \in M$,
  \begin{equation}
    \label{eqn:snorm}
    \|S[\kappa]\| \le C\|\partial_t f\|_{L^1(\bR)}.
  \end{equation}
\end{lemma}
\begin{proof} From the assumptions about $f$ and bounds
  \ref{eqn:best},  $b[\kappa](\bx_r,\cdot;\bx_s)H *
  \partial_tf \in {\cal S}(\bR)$ for all $(\bx_s,\bx_r) \in X_{sr}$,
  whence the first conclusion. Bounded extension to $L^2(X_{sr} \times
  \bR)$ follows from Young's inequality.
\end{proof}

\begin{prop}
  \label{thm:unorm}
  Under the assumptions of Lemma \ref{thm:snorm}, the least squares problem
  \ref{eqn:regs} with $\sigma > 0$ and $d\in L^2(X_{sr} \times \bR)$ has a unique solution
  $u_{\sigma}[\kappa;d] \in L^2(X_{sr} \times \bR)$ satisfying
  \begin{equation}
    \label{eqn:unorm}
    \|u_{\sigma}[\kappa;d]\| \le \frac{C}{\sigma}\|\partial_t f\|_{L^1(\bR)}\|d\|.
  \end{equation}
\end{prop}
\begin{proof}
For $\sigma > 0$, the Lax-Milgram Theorem (\cite{Yosida}, section
III.7) assures the existence of unique minimizer
$u_{\sigma}[\kappa;d]$ satisfying the  normal
equation
\begin{equation}
  \label{eqn:normal}
  (S[\kappa]^TS[\kappa] + \sigma I) u_{\sigma}[\kappa;d] = S[\kappa]^Td,
\end{equation}
and the bound follows from the inequality \ref{eqn:snorm}.
\end{proof}

\begin{prop}
  \label{thm:tl2}
  In addition to the Exponential Decay, Single Arrival, and No Return assumptions, uppose that $ Td \in L^2(X_{sr} \times \bR)$. Then
  $Tu_{\sigma}[\kappa;d]  \in L^2(X_{sr} \times \bR)$, and for each
  $(\bx_s,\bx_r) \in X_{sr}$,
  \[
    \|Tu_{\sigma}[\kappa;d](\bx_r,\cdot;\bx_s)\| \le
  \]
  \[
    \frac{C}{\sigma}\left(\left(
        \frac{\|\partial_t f\|^2_{L^1(\bR)}}{\sigma}(\|f\|_{L^1(\bR)}+\|T\partial_tf\|_{L^1(\bR)})
      \right)\|d\|\right.
  \]
  \begin{equation}
    \label{eqn:tunorm}
    \left.
    + \|\partial_t  f\|_{L^1(\bR)}(\|d (\bx_r,\cdot;\bx_s)\| +\|Td (\bx_r,\cdot;\bx_s)\|)\right).
  \end{equation}
\end{prop}

\begin{proof}
  Since $\kappa$ is the argument of $S$ and $u_{\sigma}$ throughout, drop it
  temporarily, also from the ingredients $\tau, a, $ and $b$ of the
  Hadamard decomposition. The expressions to follow pertain to
  each trace, that is, for each $(\bx_s,\bx_r) \in X_{sr}$, so also drop these
  from the notation also.
  
  Thus
  \[
    S^TS u=
    \rho^2(a\delta  +  (\check{b}\check{H})) *\delta_{-\tau}*\partial_t
    f * (a\delta  +  (bH)) *\delta_{\tau}*\partial_t f*u
  \]
  \[
    = \rho^2((a\delta  +  (\check{b}\check{H}))*(a\delta  +
    (bH))*\partial_t f *\partial_t f*.
  \]
  $T$ is also a bounded map on $ {\cal S}(\bR)$, as is
  its commutator with $S^TS$, and 
  \[
    [T,S^TS] = \rho^2((T\check{b}\check{H}) *(a\delta  +
    (bH))*+ (a\delta  +  (\check{b}\check{H}))*(TbH))*\partial_t f
    *\partial_t f*
  \]
  \[
    +2\rho^2((a\delta  +  (\check{b}\check{H}))*(a\delta  +
    (bH)))* \partial_t f  *T\partial_t f*
  \]
  \[
    =\rho^2((T\partial_t\check{b}\check{H}) *(a\delta  +
    (bH))*+ (a\delta  +  (\check{b}\check{H}))*(T\partial bH))* f
    *\partial_t f*
   \]
  \begin{equation}
    \label{eqn:tsscomm}
     +2\rho^2((a\delta  +  (\check{b}\check{H}))*(a\delta  +
    (bH)))* \partial_t f  *T\partial_t f*.
  \end{equation}
  
  Referring to the exponential decay assumption \ref{eqn:best}
  and the standing assumption $f\in C_0^{\infty}(\bR)$, 
  observe that all summands on the RHS of identity \ref{eqn:tsscomm} are in
  $L^1(\bR)$, so $[T,S^TS]$ extends to a bounded
  operator on $L^2(\bR)$, with
  \begin{equation}
    \label{eqn:tsscommnorm}
    \|[T,S^TS]\| \le C\|\partial_t f\|_{L^1(\bR)}(\|f\|_{L^1(\bR)}+\|T\partial_tf\|_{L^1(\bR)})
  \end{equation}
  
  Similarly,
  \[
    [T,S]=T(\partial_t b)H*\delta_{\tau}*f + \tau(a\delta +
    bH)*\delta_{\tau}*\partial_t f + (a\delta +
    bH)*\delta_{\tau}*T\partial_t f.
  \]
  Young's inequality shows that the RHS lies in $L^1(\bR)$, and that there exists $C$ independent of $\kappa$ so that
  $[T,S^T]$ obeys the same estimate as does $S$ (Lemma
  \ref{thm:snorm}), that is,
  \begin{equation}
    \label{eqn:tcommst}  
    \|[T,S^T]\| = \|[T,S]\| \le C\|\partial_t f\|_{L^1(\bR)}
  \end{equation}
  
  Multiplying the normal equation for $u_{\sigma}$ by $T$ and earranging,
  \[
    (S^TS+\sigma I)Tu_{\sigma} = -[T,S^TS]u_{\sigma} + [T,S^T]d + S^TTd.
  \]
  A theorem of von Neumann (\cite{Yosida}, VII.3, Theorem 2) implies
  that $S^TS+\sigma I$ has a bounded inverse, so the distribution
  $Tu_{\sigma} \in L^2(X_{sr} \times \bR)$. The inverse is bounded by
  $1/\sigma$, so the estimates \ref{eqn:snorm}, \ref{eqn:unorm},
  \ref{eqn:tcommst}, and \ref{eqn:tsscommnorm} combine to yield the
  bound \ref{eqn:tunorm}.
  
\end{proof}

The following fact will be useful later on:

\begin{prop}
  \label{thm:noisefree}
  Assume Exponential Decay, Single Arrival, and No Return. Then for
  $\log \kappa \in M$, $p \ge 1$, $F[\kappa], TF[\kappa]
  \in L^p(X_{sr} \times \bR)$, and there exist $C_0, C_1, C_2 >0$ so
  that for each $(\bx_s,\bx_r) \in X_{sr}$,
  \begin{eqnarray}
    \label{eqn:dnorms}
    C_1\|\partial_tf\|_{L^p(\bR)} - C_2\|f\|_{L^p(\bR)} \le \|F[\kappa](\bx_r,\cdot;\bx_s)\|_{L^p(\bR)} & \le & C_0\|\partial_tf\|_{L^p(\bR)}, \nonumber \\
    \|TF[\kappa](\bx_r,\cdot;\bx_s)\|_{L^p(\bR)} & \le & C_0\|\partial_tf\|_{L^p(\bR)}.
  \end{eqnarray}
\end{prop}

\begin{proof}
  A simple estimate shows that $F[\kappa]$ satisfies an exponential
  decay estimate of the form \ref{eqn:best} whence the first estimate
  follows. Since $TbH \in W^{1,1}(\bR)$,
  \[
    T\rho(a\delta + bH)*\delta_{\tau}*\partial_t f = \partial_t TbH *
    \delta_{\tau} * f + \tau \rho (a\delta +
    bH)*\delta_{\tau}*\partial_t f  + (a\delta +
    bH)*\delta_{\tau}*T\partial_t f.
  \]
  All terms involve convolution of $f$ or $\partial_t f$ with kernels
  uniformly $L^1$-bounded in $\kappa \in M$, so the second estimate follows. 
\end{proof}

{\bf Remark:} It is conventional to call data in the range of the
modeling operator ``noise-free'', and I will adopt this usage.

\subsection{Source Characteristics and Asymptotics}

Note that in the first of the two representations in display \ref{eqn:sgood}, only the
derivative $\partial_tf$ appears. For convenience, introduce
\begin{equation}
  \label{eqn:wavelet}
  w = \rho \partial_t f
\end{equation}
Then  $w \in C^{\infty}_0(\bR)$ and $\int \, w = 0$. I will also
refer to $w$ as ``the wavelet'', and view $F$ as parametrized by $w$.

Introduce asymptotics in this setting by choosing $w = \wl$
from a family $\{\wl:0 <\lambda \le 1\}$ indexed by the dimensionless parameter
$\lambda$, and defined by
\begin{equation}
  \label{eqn:wfam}
  \wl(t) = \lambda^{-1/2}w_1\left(\frac{t}{\lambda}\right).
\end{equation}
The  ``mother wavelet'' $w_1$ is assumed to have zero mean and
all wavelets $\wl$ in this family
inherit this property.
It follows from these assumptions that the Fourier transform $\fwl$ is
even and real-analytic. The corresponding source waveform $f_{\lambda}$ is
given by
\begin{equation}
  \label{eqn:flambda}
  f_{\lambda}(t) =\frac{1}{\rho} \int_{-\infty}^t\wl =
  \lambda^{1/2}f_1\left(\frac{t}{\lambda}\right).
\end{equation}
Note that for $p \ge 1$,
\begin{equation}
  \label{eqn:lambdanorms}
  \|\wl\|_{L^p(\bR)} = \lambda ^{\frac{1}{p}-\frac{1}{2}}\|w_1\|_{L^p(\bR)},\,\|f_{\lambda}\|_{L^p(\bR)}
  = \lambda ^{\frac{1}{p}+\frac{1}{2}} \|f_1\|_{L^p(\bR)}.
\end{equation}


The parameter $\lambda$ is connected to several frequency-related
characteristics of $\wl$. Suppose that $f \in L^2(\bR)$ is also square-integrable after
multiplication by $t$ (application of $T$, in the notation of the last
section - true if for instance $f$ has compact
support). The {\em pulse width} (or {\em RMS width} ) of $f$ is
\begin{equation}
  \label{eqn:pw}
  l(f) = \frac{\|Tf\|}{\|f\|}.
\end{equation}
The pulse width of $\wl$ is proportional to $\lambda$:
\begin{equation}
  \label{eqn:wpw}
  l(\wl) = \left(\frac{\int dt (t \wl(t))^2}{\int dt (\wl(t))^2}
  \right)^{1/2} = \lambda l(w_1).
\end{equation}

Note that the AWI objective defined in \ref{eqn:awidefmain} is actually
squared pulse width of the adaptive filter 
$u_{\sigma}[\kappa;d]$ respectively.

The {\em RMS frequency} of $\wl$ is
\begin{equation}
  \label{eqn:freq}
  k(\wl) = \left(\frac{\int d\omega (\omega \fwl(\omega))^2}{\int dt (w(t))^2} \right)^{1/2} = \lambda^{-1}k(w_1).
\end{equation}
Thus the {\em RMS wavelength} $1/k(\wl)$ is proportional to $\lambda$.

The Heisenberg inequality (for example, \cite{Folland:07}, p. 255) bounds the RMS
wavelength, the reciprocal of the RMS frequency, in terms of the pulse
width:
\begin{equation}
  \label{eqn:heis}
  l(\wl) k(\wl) = l(w_1)k(w_1) \ge \frac{1}{4\pi}.
\end{equation}
Equality is attained for Gaussian $w_1$, and approximated for many other
commonly used wavelets (Ricker = second derivative of Gaussian,
trapezoidal bandpass,...).

The modeling operator with wavelet $\wl$ is denoted $\Fl$
(consistently with the prior definition using $\fl$).. From the
decomposition \ref{eqn:green}, $\Fl$ is the sum of two terms:
\begin{equation}
  \label{eqn:Flamdecomp}
  \Fl = \Fla + \Flb.
\end{equation}
where
\begin{equation}
  \label{eqn:Fl0}
  \Fla[\kappa]= a[\kappa]\delta_{\tau[\kappa]}*\wl
\end{equation}
and
\[
\Flb[\kappa] = (b[\kappa]H) *\delta_{\tau[\kappa]}*\wl.
\]
Write $b_0[\kappa](\bx_s,\bx_r) = b[\kappa](\bx_r,0;\bx_s)$. Then
$\Flb$ can be rearranged as
\begin{equation}
  \label{eqn:Fl1}
  =\rho (b_0[\kappa]\delta + (\partial_t b H)) * f_{\lambda} )*\delta_{\tau[\kappa]}.
\end{equation}
From the bounds \ref{eqn:best} and \ref{eqn:lambdanorms} and Young's
inequality, there exists $C \ge 0$ so that for all $\log \kappa \in M$, $p
\ge 1$,
\begin{equation}
  \label{eqn:Frembd}
 \|\Fla[\kappa]\|_{ L^p(X_{sr} \times \bR)} \le C\lambda^{\frac{1}{p}-\frac{1}{2}},\,\|\Flb[\kappa]\|_{ L^p(X_{sr} \times \bR)} \le C\lambda^{\frac{1}{p}+\frac{1}{2}}
\end{equation}

Similarly, define
\begin{eqnarray}
  \label{eqn:sdeflam}
  \Sl[\kappa]u &=& u *\Fl[\kappa],  \nonumber \\
  \Sla[\kappa]u &=& u *\Fla[\kappa], \nonumber\\
  \Slb[\kappa]u &=& u *\Flb[\kappa], 
\end{eqnarray}
so that
\[
  \Sl[\kappa] = \Sla[\kappa] + \Slb[\kappa].
\]
From the bound \ref{eqn:Frembd} for $p=1$, there is $C \ge 0$ so that
for all $\log \kappa \in M$, $u \in L^2(X_{sr} \times \bR)$,
\begin{eqnarray}
  \label{eqn:Srembd}
  \|\Sl[\kappa]\| &\le&  C\lambda^{\frac{1}{2}},\nonumber \\
  \|\Sla[\kappa]\|& \le& C\lambda^{\frac{1}{2}}. \nonumber \\
  \|\Slb[\kappa]\| & \le& C\lambda^{\frac{3}{2}}.
\end{eqnarray}

\subsection{Linking Regularization and Wavelength}
Note that noise-free data $\dl = \Fl[\kappa]$ is bounded in $L^2(X_{sr}
\times \bR)$ uniformly in $\kappa$ and $\lambda$  (Proposition
\ref{thm:noisefree}, estimates \ref{eqn:Frembd}). 
The minimum value of the right-hand side in equation \ref{eqn:regs} is
the length of the projection of $(\dl,0)^T$ onto the range of
$(\Sl[\kappa[],\sigma I)^T$, so is less than $\|\dl^2\|$.

Suppose conversely that $\|\Sl[\kappa]\ul[\kappa;d] - \dl\| \le
\frac{1}{2}\|d\|$. Then $\|\Sl[\kappa]\ul[\kappa;d] \| \ge
\frac{1}{2}\|d\|$, so
\[
  \|d\| \le 2\|\Fl[\kappa]*\ul[\kappa;d] \| \le
  2\lambda^{\frac{1}{2}}\|\ul[\kappa;d]\|
\]
from Young's inequality, Proposition \ref{thm:noisefree}, and
inequalities \ref{eqn:lambdanorms}. From the bound \ref{eqn:unorm} this is in
turn
\[
  \le R^* \frac{\lambda}{\sigma} \|d\|
\]
whence
\[
  \sigma \le R^* \lambda.
\]
for a suitable $R^*>0$ independent of $\log \kappa \in M$ and $\lambda
\in (0,1]$. This establishes conclusion 2a in Theorem \ref{thm:main}, and
suggests that in order to expect filtered predicted data
$\Sl[\kappa]\ul[\kappa;d\|$ to remian uniformly (in $\lambda$) close to the data to be
inverted, as the heuristic justification for AWI suggests, the
regularization weight $\sigma$ must be coupled to the wavelength proxy
$\lambda$. Therefore we assume from now on that equation
\ref{eqn:lamsig} holds,  
in which $r$ is a positive constant with dimensions pressure$^2 \times 
$ time$^2$.

\subsection{The leading term and travel time}

The relation between MSWI, AWI, and travel time emerges most clearly
through a further simplification:

\begin{prop}
  \label{thm:leading}
In the definitions \ref{eqn:mswidefmain} of $J_{\rm MSWI,\sigma}$ and
\ref{eqn:awidefmain} of $J_{\rm AWI,\sigma}$, set $\sigma = r\lambda$
(equation \ref{eqn:lamsig})
and replace the data $d$ with
$\dla = \Fla[\kappa^*]$, $S[\kappa]$ with $\Sla[\kappa]$, and $\ul$
with the solution $\uzl$ of the corresponding least squares problem. That is,
\begin{eqnarray}
  \label{eqn:leading}
  \dla & = & \Fla[\kappa^*] = a[\kappa^*]\delta_{\tau[\kappa^*]}*\wl
             \nonumber \\
  \Sla[\kappa]u &=& a[\kappa]\delta_{\tau[\kappa^*]}*\wl *u\nonumber
  \\
  \uzl &=& \mbox{argmin}_{u \in L^2(X_{sr} \times \bR)} (\|\Sla[\kappa]u-\dla\|^2 + r\lambda\|u\|^2)
\end{eqnarray}
Then there exists $0 < C_* \le C^*$, $0 < C$, and $W: M \times X_{sr}
\rightarrow \bR^+$ so that for
$\kappa \in M$, $(\bx_s,\bx_r) \in X_{sr}$, $C_* \le
W[\kappa](\bx_s,\bx_r) \le C^*$, and $\uzl$ defined in display
\ref{eqn:leading} for $0 < \lambda \le 1$, 
\begin{equation}
  \label{eqn:mswittz}
 | \|T\uzl[\kappa;\dla] (\bx_r,\cdot;\bx_r)\|^2- 
  \frac{1}{\lambda} W[\kappa](\bx_s,\bx_r) (\tau[\kappa^*]({\bf
  x}_s,{\bf x}_r)-\tau[\kappa]({\bf x}_s,{\bf x}_r))^2 | \le C\lambda,
\end{equation}
and
\begin{equation}
  \label{eqn:awittz}
  \left| \frac{\|T\uzl[\kappa;\dla] (\bx_r,\cdot;\bx_r)\|^2}{
        \|\uzl[\kappa;\dla] (\bx_r,\cdot;\bx_r)\|^2} - (\tau[\kappa^*](\bx_s,\bx_r)-\tau[\kappa]({\bf x}_s,{\bf x}_r))^2\right|  \le C\lambda^2.
\end{equation}
Finally, for $\sigma = r\lambda$ with $r$ sufficiently small,
\begin{equation}
  \label{eqn:Sutz}
  \|\Sla[\kappa]\uzl[\kappa;\dzl]\| \ge \frac{1}{2}\|\dzl\|
\end{equation}
\end{prop}

\begin{proof}
  The abbreviations
  $a^*=a[\kappa^*](\bx_s,\bx_r), a=a[\kappa](\bx_s,\bx_r),
  \tau^*=\tau[\kappa^*](\bx_s,\bx_r), \tau=\tau[\kappa](\bx_s,\bx_r)$.
  will be convenient in the computations to follow. In particular, the
  source/receiver coordinates $\bx_r,\bx_r$ will be suppressed from
  the notation. Each identity and estimate to follow refers to every
  pair $(\bx_s,\bx_r) \in X_{sr}$, that is, trace-by-trace, unless
  expllicitly indiccated otherwise.,

The right-hand side of the third equation in \ref{eqn:leading} is 
\[
= \mbox{argmin}_u (\|a \wl*u*\delta_{\tau}-a^*\wl*\delta_{\tau^*}\|^2 + r\lambda\|u\|^2)
\]
which in terms of the Fourier transforms $\fwl$, $\fu$ is
\begin{equation}
  \label{eqn:leadingft}
= \mbox{argmin}_u \frac{1}{2\pi}(\|a \fwl \fu e^{-i\omega \tau}-a^*\fwl e^{-i\omega \tau^*}\|^2 + r\lambda\|\fu\|^2).
\end{equation}
The normal equation for the least squares problem \ref{eqn:leadingft} is
\[
(a^2 |\fwl|^2 +r\lambda)\hat{u} = a a^*|\fwl|^2e^{i\omega(\tau^*-\tau)}
\]
the solution of which is (omitting $\kappa, \dl$ from the arguments,
as well as $\bx_s,\bx_r$):
\begin{equation}
  \label{eqn:hatuf}
\hat{u}_{0,\lambda} = \frac{a^*}{a}\hat{g}_{\lambda,\frac{r\lambda}{a^2}} e^{i\omega(\tau^*-\tau)}
\end{equation}
where
\begin{equation}
  \label{eqn:hatgf}
\hat{g}_{\lambda,\mu} = \frac{|\fwl|^2}{|\fwl|^2 + \mu}.
\end{equation}
for $0 < \lambda \le 1, 0 < \mu$. Since $\wl$ vanishes for large $|t|$, $\fwl$ is entire analytic hence
vanishes on a set of measure zero, whence $\hat{g}_{\lambda,\frac{r\lambda}{a^2}}:$
tends to $1$ almost everywhere as $\lambda \rightarrow 0$. Hence the
inverse Fourier transform $g_{\lambda,\frac{r\lambda}{a^2}}$ tends to $\delta$ in
the sense of distributions, and
\begin{equation}
  \label{eqn:u0s}
\uzl(t) = \frac{a^*}{a}g_{\lambda,\frac{r\lambda}{a^2}}(t-(\tau^*-\tau))
\end{equation}
to a multiple of a shifted $\delta$, in fact exactly the
non-square-integrable distribution
solution of $\bar{\Sla}[\kappa]u=\dla$. However for $\lambda>0$, $\uzl[\kappa;\dla]$ is square integrable.

Since $\wl \in C_0^{\infty}(\bR)$, $\hat{g}_{\lambda,\frac{r\lambda}{a^2}}$
and $\hat{g}_{\lambda,\frac{r\lambda}{a^2}}'$ are square-integrable, whence
from equation \ref{eqn:hatuf} $(\hat{u}^0_{r\lambda})'$ is also square
integrable. Consequently $t \uzl$ is square integrable.

The following argument is a variant of one presented
by\cite{HuangSymes2015SEG}, using ideas from Song's thesis \cite[]{Song:94}. Observe that
\[
\|T\uzl[\kappa;\dla]\|^2 =  \frac{(a^*)^2}{a^2} \int dt\, t^2|g_{\lambda,\frac{r\lambda}{a^2}}(t-(\tau^*-\tau))|^2
\]
\[
=\frac{(a^*)^2}{a^2} \int dt\, (t+(\tau^*-\tau))^2|g_{\lambda,\frac{r\lambda}{a^2}}(t)|^2
\]
\begin{equation}
  \label{eqn:MSWI}
  =\frac{(a^*)^2}{a^2} \int dt\, (t^2 + 2t
  (\tau^*-\tau)+(\tau^*-\tau)^2)|g_{\lambda,\frac{r\lambda}{a^2}}(t)|^2
\end{equation}
Since $\wl$ is real, $|\fwl|$ is even, therefore
so is
$\hat{g}_{\lambda,\frac{r\lambda}{a^2}}$, therefore $g_{\lambda,\frac{r\lambda}{a^2}}$ is even, and the linear term (in t) on the
right-hand side of \ref{eqn:MSWI} vanishes, So
\begin{equation}
  \label{eqn:MSWIasym}
  \|T\uzl\|^2 =
\frac{(a^*)^2}{a^2} \int dt\, t^2|g_{\lambda,\frac{r\lambda}{a^2}}(t)|^2
+(\tau^*-\tau)^2\frac{(a^*)^2}{a^2}\|g_{\lambda,\frac{r\lambda}{a^2}}\|^2
\end{equation}
Recalling the definition \ref{eqn:pw} of pulse width, \ref{eqn:MSWIasym} can be re-written as
\begin{equation}
  \label{eqn:MSWIasym2}
 \|T\uzl\|^2 =
  \frac{(a^*)^2}{a^2}\|g_{\lambda,\frac{r\lambda}{a^2}}\|^2\left(l(g_{\lambda,\frac{r\lambda}{a^2}})^2  + (\tau^*-\tau)^2\right)
\end{equation}


From the result \ref{eqn:u0s}, for each
$(\bx_s,\bx_r) \in X_{sr}$,
\begin{equation}
  \label{eqn:u0snorm}
  \|\uzl\|^2=
  \frac{(a^*)^2}{a^2}\|g_{\lambda,\frac{r\lambda}{a^2}}\|^2.
\end{equation}
Combine equations \ref{eqn:awittz} and \ref{eqn:u0snorm} to obtain
\begin{equation}
  \label{eqn:AWIasym2}
  \frac{\|T\uzl\|^2}{\|\uzl\|^2}= \sum_{\bx_s,\bx_r} \left(l(g_{\lambda,\frac{r\lambda}{a^2}})^2  + (\tau^*-\tau)^2\right)
 \end{equation}

So the leading-term modification of each summand in the definition of $J_{\rm AWI}$ differs from mean square
travel time error by the pulse width of
$g_{\lambda,\frac{r\lambda}{a^2}}$. The analogous relation for $_{\rm
  MSWI}$ involves the norm-squared as well.

To assess these
quantities, note first that
\[
  \|g_{\lambda,\frac{r\lambda}{a^2}}\|^2 =\frac{1}{2\pi} \int d\omega\,
  \left(\frac{|\fwl|^2(\omega)}{|\fwl|^2(\omega) + \frac{r\lambda}{a^2}}\right)^2
=\frac{1}{2\pi} \int d\omega\,
  \left(\frac{|\lambda^{1/2}\hat{w}_1(\lambda
      \omega)|^2(\omega)}{|\lambda^{1/2}\hat{w}_1(\lambda
      \omega)|^2 + \frac{r\lambda}{a^2}}\right)^2
\]
\[
=\frac{1}{2\pi}\int d\omega\,
  \left(\frac{|\hat{w}_1(\lambda
      \omega)|^2}{|\hat{w}_1(\lambda
      \omega)|^2 + \frac{r}{a^2}}\right)^2
=  \frac{1}{2\pi \lambda} \int d\omega_1\,
  \left(\frac{|\hat{w}_1(\omega_1)|^2}{|\hat{w}_1(\omega_1)|^2 +
      \frac{r}{a^2}}\right)^2 
\]
\begin{equation}
  \label{eqn:g0}
  = \frac{1}{2\pi \lambda} \int d\omega\,
  \left(\frac{|a\hat{w}_1|^2}{|a\hat{w}_1|^2 + r}\right)^2 =
  \frac{1}{\lambda} \|g_{1,r}\|^2
\end{equation}
For future use, note that this identity together with \ref{eqn:u0snorm}: implies that
\begin{equation}
    \label{eqn:u0snormlam}
  \|\uzl\| = C\lambda^{-\frac{1}{2}}.
\end{equation}
Note that according to the standing convention, $C$in this identity
depends on $M$ and $w_1$ but neither on $\kappa$ nor $\lambda$.

Similarly,
\[
  \int dt\, t^2|g_{\lambda,\frac{r\lambda}{a^2}}(t)|^2 = \frac{1}{2\pi} \int d\omega
  \left|\frac{d\hat{g}_{\lambda,\frac{r\lambda}{a^2}}}{d\omega}\right|^2
  =  \frac{1}{2\pi} \int d\omega  \left(\frac{2 \mbox{Re }\bar{\hat{w}}_{\lambda}\fwl' 
    \frac{r\lambda}{a^2}}{(|\fwl|^2 + \frac{r\lambda}{a^2})^2}\right)^2
\]
\[
  =\frac{1}{2\pi}\left(\frac{r\lambda}{a^2}\right)^2 \int d\omega \frac{(2 \lambda^{1/2}\hat{w}_1(\lambda
    \omega)\lambda^{3/2}\hat{w}_1'(\lambda \omega))^2}{(\lambda|\hat{w}_1(\lambda
    \omega)|^2 + \frac{r\lambda}{a^2})^4}
\]
\[
  =\frac{1}{2\pi\lambda}\left(\frac{r\lambda}{a^2}\right)^2\int d\omega_1 \frac{(2 \lambda^2\hat{w}_1(\omega_1)\hat{w}_1'( \omega_1))^2}{(\lambda|\hat{w}_1(
    \omega_1)|^2 + \frac{r\lambda}{a^2})^4}
  =\frac{1}{2 \pi \lambda}\left(\frac{r\lambda}{a^2}\right)^2\int d\omega_1
  \frac{(2\hat{w}_1(\omega_1)\hat{w}_1'( \omega_1)
    )^2}{\left(|\hat{w}_1(\omega_1)|^2 + \frac{r}{a^2}\right)^4}
\]
\[
  =\frac{r\lambda}{2\pi a^4}\int d\omega_1
  \frac{(2\hat{w}_1(\omega_1)\hat{w}_1'( \omega_1)
    )^2}{\left(|\hat{w}_1(\omega_1)|^2 + \frac{r}{a^2}\right)^4}
\]
\begin{equation}
  \label{eqn:g1}
 = \frac{r\lambda}{2\pi} \int d\omega
  \frac{(2a\hat{w}_1a\hat{w}_1'
    )^2}{\left(|a\hat{w}_1|^2 + r\right)^4}
\end{equation}

Define
\begin{equation}
  \label{eqn:Wdef}
  W[\kappa] = \frac{(a^*)^2}{a^2}\|g_{1,r}\|^2
\end{equation}
Since $W$ depends on $\kappa$ through $a$, and on $a^*$, it satisfies
bounds of the form stated in the Proposition. 

Combining the relations \ref{eqn:g0} and \ref{eqn:g1}, obtain
\begin{equation}
  \label{eqn:pwg}
  |l(g_{\lambda,\frac{r\lambda}{a^2}})| \le C \lambda
\end{equation}

The conclusions \ref{eqn:mswittz} and \ref{eqn:awittz} follow directly
from these observations and equations \ref{eqn:MSWIasym2} and
\ref{eqn:AWIasym2}.

Similar reasoning establishes the inequality \ref{eqn:Sutz}:
\[
  \Sla[\kappa]\uzl[\kappa;\dla](\omega) = a^* \wl(\omega) \left(\frac{|\wl(\omega)|^2}{|\wl(\omega)|^2
      + \frac{\sigma}{a^2}}\right)
\]
\[
  = a^* \lambda^{\frac{1}{2}}w_1(\lambda \omega) \left(\frac{a^2|w_1(\lambda\omega)|^2}{a^2|w_1(\lambda\omega)|^2
      + \frac{\sigma}{\lambda}}\right),
\]
so
\[
  \|\Sla[\kappa]\uzl[\kappa;\dla]\|^2 = \left(\frac{a*}{a}\right)^2
  \int a^2|w_1|^2\left(\frac{a^2|w_1|^2}{a^2|w_1|^2 +
      \frac{\sigma}{\lambda}}\right)^2 \ge \frac{\lambda}{\sigma}
  (a^*)^2a^2 \int |w_1|^4 
\]
By virtue of the bound \ref{eqn:aest}, $(a^*)^2a^2 \int |w_1|^4$ is bounded below
by a positive constant independent of $\log \kappa \in M,
(\bx_s,\bx_r) \in X_{sr}$. Since $\dl$ is bounded above by a similar
constant, conclude if $r = \sigma/\lambda$ is sufficiently small, the
inequality \ref{eqn:Sutz} holds.
\end{proof}

\noindent {\bf Remark:} Inequality \ref{eqn:Sutz} is the converse
promised earlier to the observation that requiring the filtered
predicted data $\Sl[\kappa]u_{\sigma}[\kappa;\dl]$, hence the filtered
data residual, to be comparable in norm to the data $d$ implied that
$\sigma$ should be dominated by $\lambda$, as the relation
$\sigma=r\lambda$ ensures. Of course, this is an analogue, with $\Sl$
replaced by its leading approximation $\Sla$, and $\dl$ replaced by $\dla$. 

\subsection{Estimating the Remainder}

\begin{prop}
  \label{thm:remest}
  Under the hypotheses of Theorem \ref{thm:main},
  \begin{eqnarray}
    \|\ul-\uzl\| &\le& C \lambda^{\frac{1}{2}}\|\dl\|\label{eqn:diffu} \\
    \|T(\ul-\uzl)\| &\le& C \lambda^{\frac{1}{2}}\|\dl\|\label{eqn:difftu}
  \end{eqnarray}
\end{prop}

\begin{proof}
  A little algebra yields
  \begin{equation}
    \label{eqn:step1}
    ((\Sla)^T\Sla + r\lambda I) (u_{\lambda}-\uzl) = (\Slb)^T d_{\lambda} + (\Sla)^T(d_{\lambda}-d_{0,\lambda})- 
    ((\Sla)^T\Slb + (\Slb)^T\Sla + (\Slb)^T\Slb)u_{\lambda}
  \end{equation}
  Note that the bounds \ref{eqn:Srembd} imply
\[
  \|(\Sla)^T (\dl-\dzl)\| = \|(\Sla)^T(\rho \fl b^*(0)\delta +
  \rho \fl*(\partial_tb^*H)*\delta_{\tau^*}\|
\]
\[
  = \rho \|a \wl * \fl * (b^*(0)\delta - \partial_tb^*H)\|
  \le C \lambda^{\frac{3}{2}}\|\dl\|
\]
and similarly
\[
 \| (\Slb)^T\dl \| \le C \lambda^{3/2}\|\dl\|.
\]
A bit more algebra gives
\[
  (\Sla)^T\Slb + (\Slb)^T\Sla + (\Slb)^T\Slb = \rho^2 (c \delta + g) *
  \fl * \fl *,
\]
in which
\[
  c = 2a \partial t b(0) - b(0)^2)
\]
and
\[
  g = a(\partial^2_t b H + \partial^2_t \check{b}\check{H}) -
  b((0)((\partial_t \check{b}) \check{H} + (\partial_t b) H +
  ((\partial_t \check{b}) \check{H})*((\partial_t b)H)) 
\]
According to Young's inequality and the bounds \ref{eqn:best}, $g$ is integrable with a 
$L^1(X_{sr}\times \bR)$ bound independent of $\log \kappa \in
M$. Therefore from the bounds \ref{eqn:lambdanorms}, obtain
\begin{equation}
  \label{eqn:step2}
\|(\Sla)^T\Slb + (\Slb)^T\Sla + (\Slb)^T\Slb\} \le C \lambda^3
\end{equation}
From the bounds \ref{eqn:unorm} and \ref{eqn:lambdanorms},
\begin{equation}
  \label{eqn:ulnorm}
  \|\ul\| \le C \lambda^{-\frac{1}{2}}\|\dl\|
\end{equation}
Combining the preceding results in thet estimate in display
\ref{eqn:diffu}.

Substituting $\fl$ for $f$, $\dl$ for $d$ and $r\lambda$ for $\sigma$ in the estimate \ref{eqn:tunorm} and using
the bounds \ref{eqn:lambdanorms} gives
\begin{equation}
  \label{eqn:step3}
  \|T\ul\| \le C \lambda^{-\frac{1}{2}}\|\dl\|
\end{equation}
From identity \ref{eqn:step1}, obtain
\[
   ((\Sla)^T\Sla + r\lambda I) T(\ul-\uzl) = -[T,
   (\Sla)^T\Sla](\ul-\uzl) - T \Slb)^T \dl 
 \]
\begin{equation}
  \label{eqn:step4}
 + T (\Sla)^T(\dl-\dzl)- 
  T((\Sla)^T\Slb + (\Slb)^T\Sla + (\Slb)^T\Slb)\ul)
\end{equation}

To estimate the first term on the right-hand side, note that
\[
  (\Sla)^T\Sla = a^2 \wl * \wl*
\]
so
\[
  [T,  (\Sla)^T\Sla] = 2 a^2 \wl*T\wl*.
\]
Since
\begin{equation}
  \label{eqn:twl}
  \|T\wl\|_{L^1(X_{sr}\times \bR)} \le C\lambda^{\frac{3}{2}}
\end{equation}
the estimate \ref{eqn:step2} implies
\begin{equation}
  \label{step5}
  \|[T,  (\Sla)^T\Sla](\ul-\uzl)\| \le C \lambda^{\frac{5}{2}}\|\dl\|,
\end{equation}

For the second term, compute
\[
  [T, \Slb] = [T,bH*\delta_{\tau}*\wl*] = ((T+\tau)bH*\wl*
  + bH*T\wl)*\delta_{\tau}*
\]
\[
  = (((\partial_t (T+\tau)b) + \tau
  b(0)\delta)* \fl * + bH*(T\wl))*\delta_{\tau}
\]
so
\[
  \|[T,\Slb]\| \le C\lambda^{\frac{3}{2}}
\]
and so using \ref{eqn:dnorms} and \ref{eqn:lambdanorms}
\[
\|T\Slb \dl \|=  \|[T,(\Slb)^T]\dl \|+ \Slb T\dl\| \le C\lambda^{\frac{3}{2}}\|\dl\|
\]
Similarlly,
\[
\|T\Sla (\dl-\dzl) \le C \lambda^{\frac{3}{2}}\|\dl\|
\]
Finally,
\[
  [T, (\Sla)^T\Slb + (\Slb)^T\Sla + (\Slb)^T\Slb] = [T,g*\fl*\fl*] =
  Tg*\fl*\fl* + 2 g*\fl*T\fl*
\]
in the notation introduced above. Since $g$ is a linear combination of
convolutions of truncated exponentially decreasing kernels, it remains
integrable after multiplication by $t$, and
\[
\|T\fl\|_{L^1(\bR)} = C\lambda^{\frac{5}{2}}. 
\]
Since $\|T\dl\| \le C\|\dl\|$, the estimates \ref{eqn:step2} and \ref{eqn:step3} implies
\[
  \|T ((\Sla)^T\Slb + (\Slb)^T\Sla + (\Slb)^T\Slb)\ul\| \le
  C\lambda^{3}\|\ul\| \le C\lambda^{\frac{5}{2}} \|\dl\|
\]
Adding up, the right-hand side of equation \ref{eqn:step4} is bounded
by
\[
  C \lambda^{\frac{3}{2}}\|\dl\|
\]
which together with preceding inequalities establishes the second estimate \ref{eqn:difftu}.

\end{proof}

\begin{proof} {\em of Theorem \ref{thm:main}}:
Conclusions 1 and 2a have already been established. Conclusion 2b
follows from equation \ref{eqn:Sutz} in Proposition \ref{thm:leading},
combined with the first estimate in \ref{thm:remest} and the bounds
\ref{eqn:Srembd}.

  From the inequality \ref{eqn:mswittz}
  (suppressing the arguments $\bx_s,\bx_r$),
  \[
   | \|T\uzl[\kappa;\dla] \|- 
   \lambda^{-\frac{1}{2}}
   W[\kappa]^{\frac{1}{2}}|\tau[\kappa^*]-\tau[\kappa]| | \le
   C\lambda^{\frac{1}{2}}
 \]
 Combine with inequality \ref{eqn:difftu} to establish the first
 result \ref{eqn:mswitt}.

 The leading-term result \ref{eqn:awittz} implies that
 \begin{equation}
   \label{eqn:awittzsr}
   \left| \frac{\|T\uzl[\kappa;\dla] \|}{
       \|\uzl[\kappa;\dla] \|} - |\tau-\tau^*|\right|  \le C\lambda.
 \end{equation}
 A little algebra gives
 \[
   \frac{\|T\ul\|}{\|\ul\|} \le \left(\frac{\|T\uzl\|}{\|\uzl\|}+\frac{\|T\ul
     -T\uzl\|}{\|\uzl}\right)\left(1-\frac{\|\ul-\uzl\|}{\|\uzl\|}\right)^{-1}
 \]
 \[
   \le \left(\frac{\|T\uzl\|}{\|\uzl\|} + C\lambda \|\dl\|\right) (1-C\lambda)
 \]
 Using the result \ref{eqn:u0snormlam}, and assuming $\lambda \le
 1/C$, find
 \[
 \frac{\|T\ul\|}{\|\ul\|}-\frac{\|T\uzl\|}{\|\uzl\|} \le C\lambda
 \left(1+2 \frac{\|T\uzl\|}{\|\uzl\|} \right) + 2 C^2 \lambda^2
\]
\[
  \le \le C\lambda  (1+2|\tau - \tau^*| + 2 (C^2+1) \lambda^2
\]
Since the maximum difference of traveltimes is a function of $M$, this
can be written as
\[
  \le C\lambda
\]
Similarly,
\[
  \frac{\|T\uzl\|}{\|\uzl\|}-\frac{\|T\ul\|}{\|\ul\|} \le C\lambda
\]
Combine these inequalities with the leading-term result
\ref{eqn:awittzsr} to yield the estimate \ref{eqn:awitt}.

\end{proof}

\section{Relation to Penalty Formulation}
Penalty functions provide a means to convert naively posed
inverse problems into well-posed problems, with solutions that exist
and depend stably on proper data. The ``pre-whitening'' construction
explained in the last section is a classic example, adding a scaled
square-integral of the adaptive filter to the residual square-integral to
produce an objective with a square-integrable (finite-energy)
minimizer. This penalty construction if often called {\em Tihonov
  regularization}. The MSWI and AWI objective functions have a deeper
relation to a penalty construction, however: they are
suitably defined limits of penalty functions.

This relation is abstract, and does not only pertain to the problems
discussed in this paper. To derive this abstract relation, suppose for
the moment that $S:U \rightarrow D$ is a bounded and coercive map from a
domain Hilbert space $U$ to a range Hilbert space $D$. That is, for
real numbers $0 < C_* \le C^*$,
\[
  C_*\|u\|_U \le \|Su\|_D \le C^*\|u\|
\]
Suppose that $d \in D$ (``data''), and $T$ is another bounded operator $U \rightarrow V$, $V$
another Hilbert space. Define
\begin{equation}
  \label{eqn:eq1}
  J_{\alpha}(u) = \frac{1}{2}(\|Su-d\|_D + \alpha^2\|Tu\|_V^2).
\end{equation}
and
\begin{equation}
  \label{eqn:eq0}
  \tJa= \min_u J_{\alpha}(u) = J_{\alpha}(u_{\alpha})
\end{equation}
Since $S$ is coercive, the minimum is well-defined for any
$\alpha \ge 0$, and the minimizer $\ua \in U$ satisfies the normal
equation
\begin{equation}
  \label{eqn:eq0n}
  (S^TS + \alpha^2 T^TT) \ua = S^Td
\end{equation}

The claim to be established is
\begin{prop}
  \label{thm:alphalim}
  Suppose that $U$, $V$ and $D$ are Hilbert spaces,
  $S:U\rightarrow D$ is bounded and coercive, $T:U\rightarrow V$ is
  bounded, $\ua$ is the solution of
  the normal equation \ref{eqn:eq0n}, and $J_{\alpha}$,
  and $\tJa$ are defined by \ref{eqn:eq1} and \ref{eqn:eq0}
  respecively for $\alpha \ge 0$. Then
\begin{equation}
  \label{eqn:eq4}
  \lim_{\alpha \rightarrow 0} \frac{1}{\alpha^2}  (\tJa-\tJz).= \frac{1}{2}\|T\uz\|^2
\end{equation}
\end{prop}

\begin{proof}
Use the normal equation to write
\[
  \ua = (S^TS + \alpha^2 T^TT)^{-1}S^Td = (S^TS + \alpha^2 T^TT)^{-1}S^TS \uz
\]
\begin{equation}
  \label{eqn:eq2}
  = \uz - \alpha^2 (S^TS + \alpha^2 T^TT)^{-1}T^TT\uz = \uz-\alpha^2 \va
\end{equation}
Note that $\va = (S^TS + \alpha^2 T^TT)^{-1}T^TT\uz$ is uniformly bounded in $\alpha \ge 0$.

Substitute the RHS of equation \ref{eqn:eq2} into the definition \ref{eqn:eq1} to obtain
\begin{equation}
  \label{eqn:eq3}
  \tJa = \frac{1}{2}(\|S\uz-d\|_D^2 - 2 \alpha^2\langle S\uz-d, S\va\rangle_D + \alpha^2\|T\uz\|_V^2 + O(\alpha^4))
\end{equation}
The second term on the RHS of equation \ref{eqn:eq3} vanishes thanks to the normal equation, and the first term is precisely $\tJz$. The conclusion \ref{eqn:eq4} follows.
\end{proof}

To apply the assertion \ref{eqn:eq4} to the MSWI objective, make the
following choices, using the notation developed in the last section.
In particular, $\sigma$ and $\lambda$ are related throughout by
$\sigma=r\lambda$, where $r$ is a fixed positive constant with unit
pressure$^2 \times$time$^2$.
For each $\kappa \in M$, $\lambda \in (0,1]$, $S$ is the regularized
modeling operator $[\Sl,\sqrt{r\lambda} I]^T$, with domain is
$U=L^2(X_{sr}\times\bR)$ and range $D =  L^2(X_{sr} \times \bR) \oplus
L^2(\bR \times \bR)$. The role of $d$ is played by the
augmented data $(d,0)^T$, $d \in D$. $T$ is a modified multiply-by-$t$ operator
$\Te$, defined by 
\begin{equation}
  \label{eqn:teps}
  \Te u(t) = \frac{t}{\sqrt{1+\epsilon^2t^2}}.
\end{equation}
with domain and range $U ( = V)$. $\Te$ is bounded for $\epsilon>0$,
and $=T$ on any natural dense domain as in the last section for $\epsilon=0$. The final result below will
concern the limit $\epsilon \rightarrow 0$.

Then $\Ja$ becomes $\Jase[\kappa,u;d]$:
\begin{equation}
  \label{eqn:mswipen}
   \Jase[\kappa,u;d] = \frac{1}{2}\left(\|\Sl[\kappa]u-d\|^2 +
   r\lambda\|u\|^2 + \alpha^2\|\Te u\|^2\right)
 \end{equation}
The minimizer $\ua$ over $u$ of $\Jase[\kappa,u;d]$ is
$\uase[\kappa;d]$, and $\tJa$ becomes
\begin{equation}
  \label{eqn:mswipenred}
  \tJase[\kappa;d] =
  \Jase[\kappa,\uase[\kappa;d];d].
\end{equation}

Define
\begin{equation}
  \label{eqn:tj0def}
  \tilde{J}_{\sigma}[\kappa;d] = \|S[\kappa]\ul[\kappa;d] - d\|^2 + \sigma\|\ul[\kappa;d\|^2.
\end{equation}
in which $\ul$ is the solution of the regularized least squares
problem \ref{eqn:reg}. 

\begin{prop}
  \label{thm:MSWIepsalphalim}
Suppose that $\kappa^* \in M$ and $\dl=\Fl[\kappa^*]$.
\begin{equation}
  \label{eqn:eq7}
  \lim_{\epsilon \rightarrow 0}\lim_{\alpha \rightarrow 0}
  \frac{1}{\alpha^2}\left(\tJase[\kappa;\dl]-\tilde{J}_{\sigma}[\kappa;d]\right)=
 \|T\ul\|^2 =  J_{\rm MSWI,\sigma}[\kappa;d]
\end{equation}
in which $\ul$ is the solution of the least squares problem
\ref{eqn:reg} with $S=\Sl$, $d=\dl$, and $J_{\rm MSWI,\sigma}$ is as
defined in equation \ref{eqn:mswidefmain}.
\end{prop}

\begin{proof}
  It follows immediately from Proposition \ref{thm:alphalim} that
  \[
    \lim_{\alpha \rightarrow 0}
    \frac{1}{\alpha^2}\left(\tJase[\kappa;\dl]-\tilde{J}_{\rm MSWI,0,\sigma,\epsilon}[\kappa;\dl]\right)=
    \|\Te u_{\rm MSWI,0,\sigma,\epsilon}\|^2
  \]
  Comparing the normal equation \ref{eqn:eq2} for $\alpha=0$
  to the definition \ref{eqn:reg}, one sees that  $u_{\rm
    MSWI,0,\sigma,\epsilon} = u_{\rm MSWI,0,\sigma,0} = \ul$ is
  independent of $\epsilon$. From
  Propositions \ref{thm:tl2} and \ref{thm:noisefree}, $T\ul \in
  L^2(X_{sr} \times \bR)$. For $\delta>0$, choose $t_{\delta} > 0$ so
  that
  \[
    \int_{|t| \ge t_{\delta}} \,dt\, |\Te \ul|^2 \le  \int_{|t| \ge
      t_{\delta}} \,dt\, |T \ul|^2 \le \delta^2
  \]
  and $\epsilon_{\delta}>0$ so that
  \[
    \int_{|t| \le t_{\delta}} \,dt\, |(T-\Te) \ul|^2 \le \delta^2
  \]
  for $0<\epsilon<\epsilon_{\delta}$. Adding up, obtain $\|T\ul
  -\Te\ul\|^2 \le 5\delta^2$. Since $\delta>0$ was arbitrary, it
  follos that $\|(\Te -T)\ul\| \rightarrow 0$ as $\epsilon \rightarrow
  0$..
\end{proof}

\begin{cor}
  \label{thm:cor2}
  With the hypotheses and notation of Theorem \ref{thm:main},
  \[
    \lim_{\lambda \rightarrow 0}\lim_{\epsilon \rightarrow 0}\lim_{\alpha \rightarrow 0}
  \frac{\lambda}{\alpha^2}\left(\tilde{J}_{\rm MSWI,\alpha,r\lambda,\epsilon}[\kappa;\dl]-\tilde{J}_{r\lambda}[\kappa;\dl]\right)=
  \]
  \begin{equation}
    \label{eqn:main1}
    = \sum_{(\bx_s,\bx_r) \in X_{sr}}
    M[m,m^*,w_1]^2(\bx_s,\bx_r) |\tau[\kappa](\bx_s,\bx_r) -
    \tau[\kappa^*](\bx_s,\bx_r)|^2.
  \end{equation}
\end{cor}

To derive the analogous relation for AWI, define a family of norms on
$ L^2(X_{sr} \times \bR)$, depending on $\log \kappa \in M, d \in
L^2(X_{sr} \times \bR)$, and $\sigma >0$:
\begin{equation}
  \label{eqn:awinormdef}
  \|u\|^2_{\kappa,d,\sigma} = \sum_{(\bx_s,\bx_r) \in X_{sr}}\frac{\|u(\bx_r,\cdot;\bx_s)\|^2}{\|\ul[\kappa;d](\bx_r,\cdot;\bx_s)\|^2}
\end{equation}
Define $V_{\kappa,d,\sigma} = L^2(X_{sr} \times \bR)$, equipped with
the norm $\|\cdot\|_{\kappa,d,\sigma}$, and define $\Tekds$ to be
$\Te$, with domain $U$ and range $V_{\kappa,d,\sigma}$.

The AWI penalty objective is then
\begin{equation}
  \label{eqn:awipen}
  \JAase[\kappa,u;d] = \|S[\kappa]u-d\|^2 + \sigma\|u\|^2 + \alpha^2
  \|\Tekds u\|_{\kappa,d,\sigma}^2
\end{equation}
As before, define $\uaase[\kappa;d]$ to be the minimizer of
$\JAase[\kappa,u;d]$ over $u \in U$, and
\begin{equation}
  \label{eqn:awipenred}
  \tJAase[\kappa;d] =\JAase[\kappa,\uaase[\kappa;d];d].
\end{equation}
An argument precisely analogous to the proof of Corollary
\ref{thm:cor2} leads to the conclusion of
\begin{cor}
  \label{thm:cor3}
  With the hypotheses and notations of Theorem \ref{thm:main}, and the
  definitions \ref{eqn:tj0def} and \ref{eqn:awipenred} with $S=\Sl$,
  $d=\dl$, and $\sigma=r\lambda$,
  \[
  \lim_{\lambda \rightarrow 0}\lim_{\epsilon \rightarrow 0}\lim_{\alpha \rightarrow 0}
  \frac{1}{\alpha^2}\left(\tilde{J}_{\rm AWI,\alpha,r\lambda,\epsilon}[\kappa;\dl]-\tilde{J}_{r\lambda}[\kappa;\dl]\right)=
  \]
  \begin{equation}
    \label{eqn:cor3}
    = \sum_{(\bx_s,\bx_r) \in X_{sr}} |\tau[\kappa](\bx_s,\bx_r) - \tau[\kappa^*](\bx_s,\bx_r)|^2.
  \end{equation}
\end{cor}

\noindent {\bf Remark:} From inequalities \ref{eqn:ulnorm},
\ref{eqn:Srembd}, $\tilde{J}_{r\lambda}[\kappa;\dl] \le C$. That
suggests that the contribution of the second term on the
left-hand-side of the limit \ref{eqn:cor1mswi} is insignificant for small
$\lambda$, whereas it is comparable to the first term in
\ref{eqn:cor3}.


\section{Multiple Arrivals and Cycle-Skipping}
The geometric asymptotics approximation \ref{eqn:green} holds only
in the ``single arrival'' case, that is, when sources and receivers
are close enough together relative to the length scale of significant
variation in the wave velocity $((\kappa\beta)^{1/2})$ that a unique ray of
geometric acoustics connects each source-receiver pair. At larger
distances, multiple raypaths typically exist, corresponding to
multiple idenfiable wavefront arrivals \cite[]{Whi:82}. Accordingly, in such settings
the analysis presented in the last section no longer applies, and the
relation established there between AWI and travel-time inversion is
case into doubt.

A simple calculation reveals the phenomenon discussed in all of these
works. I will argue formally, via the leading term in the Hadamard
approximation of the Green's function \ref{eqn:green}, generalized to
multiple arrival times. The analogue of the geometric asymptotics for a source-receiver
pair with connected by multiple rays is (once again suppressing the
source and receiver coordinates from the notation):
\begin{equation}
  \label{eqn:multi}
  F[\kappa](t) \approx \sum_{i \ge 0} a_i H^{p_i}w(t-\tau_i).
\end{equation}
The amplitudes $\{a_i\}$ and arrival times $\{\tau_i\}$ depend on $\kappa$, as
before, through ray-tracing. $H$ is the Hilbert transform, and $p_i$
counts the number of times the ray for arrival time $\tau_i$ has
touched a caustic (counted with multiplicity, in an appropriate
sense). The earliest arrival has not touched a caustic. If the
indexing is organized so that $\tau_0 < \tau_i$ for $i>0$, then
$p_0=0$. See \cite{Friedlander:75} for details.

The right-hand side of equation \ref{eqn:multi} can be expressed as
\[
  = (\delta + T) * a_0w(t-\tau_0),
\]
in which 
\[
  T = \sum_{i \ge 1}
  \frac{a_i}{a_0}H^{p_i}\delta(t-(\tau_i-\tau_0))
\]
Thus
\begin{equation}
  \label{eqn:unwrap}
  (\delta +T_1) * F[\kappa,w] \approx  a_0w(t-\tau_0)
\end{equation}
where $\delta + T_1$ is the convolution inverse of $\delta + T$,
defined by 
\[
  T_1 = \sum_{n=1}^{\infty} (-1)^n T^n.
\]
$T^n$ should be understood as the $n$th convolution power in this
expression, that is, $T*T*...*T$, $n$ times. $T_1$ is well-defined as a distribution.

Once again, it is possible to write the adaptive filter u for which $S_{\lambda}[\kappa]u = F[\kappa]*u = d =
F[\kappa^*]$ by inspection. Use the notations $a^*_i, t^*_i, p^*_i$ for
the analogous quantities computed with the ``true'' model $\kappa^*$, and
$T^*, T_1^*$ defined as above but with $a^*_i$ in place of $a_i$ and so
on:
\[
  d = F[\kappa_*]= S[\kappa] u
\]
provided that we choose
\begin{equation}
  \label{eqn:notl2again}
u(t) \approx (\delta + T^*)*(\delta +
T_1)*\left(\frac{a^*_0}{a_0}\delta(t - (\tau^*_0-\tau_0))\right)
\end{equation}
Recalling the ordering of arrival times, 
\[
  T^* = \frac{a^*_1}{a^*_0}H^{p_1^*} \delta(t-(\tau^*_1-\tau^*_0)) + ...
\]
\[
  T_1 = -\frac{a_1}{a_0}H^{p_1} \delta(t-(\tau_1-\tau_0)) + ...
\]
where the elided terms all involve $\delta$s located at later times.
So
\[
  (\delta + T^*)*(\delta + T_1)(t) = \delta(t)  +  \frac{a^*_1}{a^*_0}
  H^{p_1^*}\delta(t-(\tau^*_1-\tau^*_0)) -\frac{a_1}{a_0}H^p_1
  \delta(t-(\tau_1-\tau_0)) + ...
\]
and
\begin{equation}
  \label{eqn:notl2explicit}
  u(t) = \frac{a^*_0}{a_0}\left(\delta(t - (\tau^*_0-\tau_0)) +
    \frac{a^*_1}{a^*_0}H^{p_1^*}\delta(t-(\tau^*_1-\tau_0)) -
    \frac{a_1}{a_0}H^{p_1}\delta(t-(\tau_1-\tau_0^*)) + ...\right)
\end{equation}
The elided terms involve only travel time differences other than the
ones appearing here, so nothing cancels.

Observe that the travel-time differences $\tau^*_1-\tau_0$,
$\tau_1-\tau_0^*$ do not tend to zero, though they do tend to the same
limit as $\kappa \rightarrow \kappa^*$ (as do the amplitude quotients), and the
Hilbert transform powers $p_1, p_1^*$ are the same as soon as $\kappa$ is
sufficiently close to $\kappa^*$.

Of courese, the adaptive filter just computed is not square
integrable, and must be regularized for use in a square-integral
construction like that of AWI or MSWI. That amounts essentially to
replacing the $\delta$s in the above formulae by a square-integrable
approximate Dirac function, of roughly a wavelength in width. In the
preceding section, this approximate Dirac was $g$, so use the same
notation here. The regularized version of $u(t)$ (solution of the
least squares problem \ref{eqn:reg}) contains summands with multiples of $g(t-(\tau^*_1-\tau_0))$
and $g(t-(\tau_1-\tau^*_0))$. These are approximately scaled by
$\tau^*_1-\tau_0$ and $\tau_1-\tau^*_0$ respectively in the
definitions \ref{eqn:mswidefmain} and \ref{eqn:awidefmain} of the MSWI and AWI
objectives. For other choices of objective, they contribute similarly
to the total. They will eventually begin to cancel as $\kappa$ approaches
$\kappa^*$, however this cancellation happens only when the times are
within the width of the approximate 
Dirac waveform $g$ - that is, within a wavelength.

\section{Conclusion}


This paper has explained the behaviour of AWI (and MSWI) applied to
acoustic transmission through smoothly varying material models. For
isotropic point source/receiver data exhibiting a single smooth
wavefront, these formulations produce objective functions that converge
to versions of travel time mean-square error, in the limit of hgh
data frequency. 


Some of the assumptions leading to this conclusion 
could certainly be relaxed. For instance, neither absence of boundary effects
nor exponential decay of local energy are essential. Instead
the finite time duration of actual measured data, together with more
elaborate ray geometry assumptions, can replace the infinite time data setting
posed here,. The conclusions of this paper likely hold in some form
even for non-smooth material parameter variation, provided that the data for such models is
dominated by ballistic transmitted waves. Similar conclusions should
hold for multipole point sources and
receivers, modeling non-isotropic experimental devices that are small on the
wavelength scale.

The exception is the single-arrival hypothesis: it is
essential, as without it, the link between the travel time and
extended waveform inversion is lost.
A number of authors have illustrated this disconnect via numerical
example. The first such examples for a version of MSWI appeared in
\cite{Symes:94c}. Crosswell seismic surveys often generate multiple
arrivals near low-velocity sedimentary layers that are natural targets
of seismic exploration and act as wave guides. \cite{Plessix:00}
successfully applied a version of MSWI to crosswell waveform inversion
by removing traces with close source and receiver depths, thus leaving
only data that tended to be dominated by single arrivals.
\cite{HuangSymes:Geo17} and \cite{SymesChenMinkoff:24} illustrate the
failure of the version of MSWI reported here to avoid cycle-skipping
in inversion of strongly refracting models, using frequency-domain
implementation of the basic definition \ref{eqn:mswidefmain} and
time-domain implementation of the reduced penalty function
\ref{eqn:mswipenred}, respectively. \cite{Yongetal:GJI23} show that
precisely the same shortcoming afflicts AWI.

Overcoming this limitation requires substantial change in the
extension of basic acoustic modeling underlying AWI. For example,
\cite{Yongetal:GJI23} have proposed localizing the definition of the
AWI objective to short time intervals, within which multiple wavefront
arrrivals are unlikely to appear, and illustrated the effectiveness of
this approach by numerical example. \cite{Symes:23} constructs an
extension based on surface sources, as opposed to the trace-by-trace
extension which is the basis of MSWI and AWI. This surface source
extension separates events by infinitesimal slope, information not
available in the trace-by-trace construction..

The analysis presented here has nothing to say about AWI applied to
reflected wave data. \cite{Warner:16} and others have
suggested that AWI is effective in that context. If so, the underlying
mechanism remains to be explained.

Linear elasticity is a much better model for mechanical waves
propagating in some human tissue (notably bone) and the solid earth
than is acoustodynamics. Least-squares FWI of elastodynamic
transmission data also tends to stagnate due to cycle-skipping, and it
is tempting to think that extension approaches such as AWI should
ameliorate that obstacle. It will be necessary extract wave speed
informaiton in the presence of potential interference of cross-talk
between wave modes, similar to cross-talk between multiple arrivals in
the acoustic case.

Finally, virtually all materials supporting wave propagation also
attenuate energy through various physical mechanisms, which act per
cycle and therefore suppress high frequency components. Analysis such
as that presented here, relying on the infinite-frequency limit
(geometric optics), is {\em a priori} impossible for such attenuative
wave physics without fundamental
reformulation..

\bibliographystyle{seg}
\bibliography{../../bib/masterref}

\end{document}